\documentstyle[12pt]{article}
\begin{document}

\title{Variational inequalities}

\author{Nikolaos E. Sofronidis\footnote{$A \Sigma MA:$ 130/2543/94}}

\date{\footnotesize Department of Economics, University of Ioannina, Ioannina 45110, Greece.
(nsofron@otenet.gr, nsofron@cc.uoi.gr)}

\maketitle

\begin{abstract}
If $- \infty < \alpha < \beta < \infty $ and $f \in C^{3} \left( [
\alpha , \beta ] \times {\bf R}^{2} , {\bf R} \right) $ is
bounded, while $y \in C^{2} \left( [ \alpha , \beta ] , {\bf R}
\right) $ solves the typical one-dimensional problem of the
calculus of variations to minimize the function $$F \left( y
\right) = \int_{ \alpha }^{ \beta }f \left( x, y(x), y'(x) \right)
dx,$$ then for any ${\phi } \in C^{2} \left( [ \alpha , \beta ] ,
{\bf R} \right) $ for which ${\phi }^{(k)}( \alpha ) = {\phi
}^{(k)}( \beta ) = 0$ for every $k \in \{ 0, 1, 2 \} $, we prove
that $\int_{\alpha }^{\beta } \left( \frac{ {\partial }^{2}f }{
\partial y^{2} } {\phi }^{2} - \frac{ {\partial }^{3}f }{
\partial y^{2} \partial y' } 2 {\phi }^{3} \right) dx$
$\geq \int_{\alpha }^{\beta } \left( \frac{ {\partial }^{2}f }{
\partial y \partial y' } 2 \phi \phi ' + \frac{ {\partial
}^{3}f }{ \partial y {\partial y'}^{2} } 2 {\phi }^{2} \phi ' +
\frac{ {\partial }^{2}f }{ {\partial y'}^{2} } \phi \phi '' +
\frac{ {\partial }^{3}f }{ \partial y {\partial y'}^{2} } \phi '
{\phi }^{2} + \frac{ {\partial }^{3}f }{ {\partial y'}^{3} } \phi
{\phi '}^{2} \right) dx$, so either the above are variational
inequalities of motion or the Lagrangian of motion is not $C^{3}$.
\end{abstract}

\section*{\footnotesize{{\bf Mathematics Subject Classification:} 49J40, 70H03, 70H30.}}

{\bf 1. Definition.} If $- \infty < \alpha < \beta < \infty $ and
$f : [ \alpha , \beta ] \times {\bf R}^{2} \rightarrow {\bf R}$ is
any bounded continuous function, then the typical one-dimensional
problem of the calculus of variations is to minimize the function
$F$, which is defined by the relation $$F(y) =
\int_{\alpha}^{\beta } f \left( x,y(x),y'(x) \right) dx, \eqno
(1)$$ where $y : [ \alpha , \beta ] \rightarrow {\bf R}$ ranges
over a suitably chosen class of functions.
\\ \rm \\
The following two propositions are well-known. See Appendix D on
pages 151-152 of [1].
\\ \rm \\
{\bf 2. Proposition.} If $- \infty < \alpha < \beta < \infty $ and
$r$ is any positive integer, while $y : \left[ \alpha , \beta
\right] \rightarrow {\bf R}$ is any continuous function such that
$$\int_{\alpha }^{\beta }y(x)\eta (x)dx=0 \eqno (2)$$ for every $\eta \in
C^{r} \left( [ \alpha , \beta ] , {\bf R} \right) $ for which
$${\eta }^{(k)}( \alpha ) = {\eta }^{(k)}( \beta ) = 0 \eqno (3)$$
for every $k \in \{ 0, 1, ..., r \} $, then $$y=0 \eqno (4)$$ on
$[ \alpha , \beta ]$.
\\ \rm \\
{\bf 3. Proposition.} If $- \infty < \alpha < \beta < \infty $ and
$y : \left[ \alpha , \beta \right] \rightarrow {\bf R}$ is any
continuously differentiable function that solves the typical
one-dimensional problem of the calculus of variations to minimize
the function $$F \left( y \right) = \int_{ \alpha }^{ \beta }f
\left( x, y(x), y'(x) \right) dx, \eqno (5)$$ where $f : [ \alpha
, \beta ] \times {\bf R}^{2} \rightarrow {\bf R}$ is $C^{2}$, then
$$\frac{ \partial f}{ \partial y } - \frac{d}{dx} \left( \frac{
\partial f}{ \partial y' } \right) = 0. \eqno (6)$$
\\
{\bf 4. Definition.} Keeping the notation and the assumptions as
in the previous proposition, given any ${\phi } \in C^{2} \left( [
\alpha , \beta ] , {\bf R} \right) $ for which $${\phi }^{(k)}(
\alpha ) = {\phi }^{(k)}( \beta ) = 0 \eqno (7)$$ for every $k \in
\{ 0, 1, 2 \} $, it is not difficult to see that the function $$I
: {\bf R} \ni t \mapsto F \left( y + t \phi \right) \in {\bf R}
\eqno (8)$$ is twice differentiable and attains its minimum at $t
= 0$, so, by virtue of Proposition 2 of Section 17.1.3 on page 409
of [2], as in Section 10.7 on page 458 of [3] or in Section 1.3 on
page 16 of [4], one obtains that
\begin{enumerate}
\item[ ]
$I''(t) = \int_{\alpha }^{\beta } \left( \frac{ {\partial }^{2}f
}{ \partial y^{2} } \left( x, y(x) + t \phi (x), y'(x) + t \phi '
(x) \right) \phi (x)^{2} \right.$
\item[ ]
$+ \ \frac{ {\partial }^{2}f }{ \partial y \partial y' } \left( x,
y(x) + t \phi (x), y'(x) + t \phi ' (x) \right) 2 \phi (x) \phi
'(x)$
\item[ ]
$\left. + \ \frac{ {\partial }^{2}f }{ {\partial y'}^{2} } \left(
x, y(x) + t \phi (x), y'(x) + t \phi ' (x) \right) \phi ' (x)^{2}
\right) dx$. \hfill (9)
\end{enumerate}

\noindent Our first purpose in this article is to prove the
following.
\\ \rm \\
{\bf 5. Proposition.} Keeping the notation as in the previous
definition, if $f$ is $C^{3}$ and $y$ is $C^{2}$, then
\begin{enumerate}
\item[ ]
$I''(0) = \int_{\alpha }^{\beta } \left( \frac{ {\partial }^{2}f
}{ \partial y^{2} } {\phi }^{2} - \frac{ {\partial }^{3}f }{
\partial y^{2} \partial y' } 2 {\phi }^{3} \right) dx
- \int_{\alpha }^{\beta } \left( \frac{ {\partial }^{2}f }{
\partial y \partial y' } 2 \phi \phi ' \right.$
\item[ ]
$\left. + \ \frac{ {\partial }^{3}f }{ \partial y {\partial
y'}^{2} } 2 {\phi }^{2} \phi ' + \frac{ {\partial }^{2}f }{
{\partial y'}^{2} } \phi \phi '' + \frac{ {\partial }^{3}f }{
\partial y {\partial y'}^{2} } \phi ' {\phi }^{2} + \frac{
{\partial }^{3}f }{ {\partial y'}^{3} } \phi {\phi '}^{2} \right)
dx.$ \hfill (10)
\end{enumerate}

\noindent {\bf Proof.} By virtue of Proposition 2 of Section
17.1.3 on page 409 of [2], one obtains that
\begin{enumerate}
\item[ ]
$I''(0) = \int_{\alpha }^{\beta } \left( \frac{ {\partial }^{2}f
}{ \partial y^{2} } {\phi }^{2} + \frac{ {\partial }^{2}f }{
\partial y \partial y' } 2 \phi \phi ' + \frac{ {\partial }^{2}f }{ {\partial y'}^{2} } {\phi '}^{2} \right)
dx$ \hfill (11)
\item[ ]
$= \ \int_{\alpha }^{\beta } \frac{ {\partial }^{2}f }{ \partial
y^{2} } {\phi }^{2} + \int_{\alpha }^{\beta } \left( \frac{
{\partial }^{2}f }{
\partial y \partial y' } 2 \phi + \frac{ {\partial }^{2}f }{ {\partial y'}^{2} } \phi '
\right) \phi ' dx$ \hfill (12)
\item[ ]
$= \ \int_{\alpha }^{\beta } \frac{ {\partial }^{2}f }{ \partial
y^{2} } {\phi }^{2} + \left[ \left( \frac{ {\partial }^{2}f }{
\partial y \partial y' } 2 \phi + \frac{ {\partial }^{2}f }{ {\partial y'}^{2} } \phi '
\right) \phi \right] _{\alpha }^{\beta } - \int_{\alpha }^{\beta }
\left( \frac{ {\partial }^{2}f }{
\partial y \partial y' } 2 \phi + \frac{ {\partial }^{2}f }{ {\partial y'}^{2} } \phi '
\right) ' \phi dx$ \hfill (13)
\item[ ]
$= \ \int_{\alpha }^{\beta } \frac{ {\partial }^{2}f }{ \partial
y^{2} } {\phi }^{2} - \int_{\alpha }^{\beta } \left( \frac{
{\partial }^{2}f }{
\partial y \partial y' } 2 \phi + \frac{ {\partial }^{2}f }{ {\partial y'}^{2} } \phi '
\right) ' \phi dx$ \hfill (14)
\item[ ]
$= \ \int_{\alpha }^{\beta } \frac{ {\partial }^{2}f }{ \partial
y^{2} } {\phi }^{2} - \int_{\alpha }^{\beta } \left( \frac{
{\partial }^{2}f }{
\partial y \partial y' } 2 \phi ' + \left( \frac{ {\partial }^{3}f
}{ \partial x \partial y \partial y' } \cdot 0 + \frac{ {\partial
}^{3}f }{ {\partial y}^{2} \partial y' } \phi + \frac{ {\partial
}^{3}f }{ \partial y {\partial y'}^{2} } \phi ' \right) 2\phi
\right.$
\item[ ]
$\left. + \ \frac{ {\partial }^{2}f }{ {\partial y'}^{2} } \phi ''
+ \left( \frac{ {\partial }^{3}f }{ \partial x {\partial y'}^{2} }
\cdot 0 + \frac{ {\partial }^{3}f }{ \partial y {\partial y'}^{2}
} \phi + \frac{ {\partial }^{3}f }{ {\partial y'}^{3} } \phi '
\right) \phi ' \right) \phi dx$ \hfill (15)
\item[ ]
$= \ \int_{\alpha }^{\beta } \frac{ {\partial }^{2}f }{ \partial
y^{2} } {\phi }^{2} - \int_{\alpha }^{\beta } \left( \frac{
{\partial }^{2}f }{
\partial y \partial y' } 2 \phi ' + \left( \frac{ {\partial
}^{3}f }{ {\partial y}^{2} \partial y' } \phi + \frac{ {\partial
}^{3}f }{ \partial y {\partial y'}^{2} } \phi ' \right) 2\phi
\right.$
\item[ ]
$\left. + \ \frac{ {\partial }^{2}f }{ {\partial y'}^{2} } \phi ''
+ \left( \frac{ {\partial }^{3}f }{ \partial y {\partial y'}^{2} }
\phi + \frac{ {\partial }^{3}f }{ {\partial y'}^{3} } \phi '
\right) \phi ' \right) \phi dx$ \hfill (16)
\item[ ]
$= \ \int_{\alpha }^{\beta } \left( \frac{ {\partial }^{2}f }{
\partial y^{2} } {\phi }^{2} - \frac{ {\partial }^{3}f }{
\partial y^{2} \partial y' } 2 {\phi }^{3} \right) dx
- \int_{\alpha }^{\beta } \left( \frac{ {\partial }^{2}f }{
\partial y \partial y' } 2 \phi \phi ' \right.$
\item[ ]
$\left. + \ \frac{ {\partial }^{3}f }{ \partial y {\partial
y'}^{2} } 2 {\phi }^{2} \phi ' + \frac{ {\partial }^{2}f }{
{\partial y'}^{2} } \phi \phi '' + \frac{ {\partial }^{3}f }{
\partial y {\partial y'}^{2} } \phi ' {\phi }^{2} + \frac{
{\partial }^{3}f }{ {\partial y'}^{3} } \phi {\phi '}^{2} \right)
dx.$ \hfill (17)
\end{enumerate}

\noindent Our second purpose in this article is to prove the
following.
\\ \rm \\
{\bf 6. Proposition.} Keeping the notation as in the previous
proposition, for any such $\phi $, we have that
\begin{enumerate}
\item[ ]
$\int_{\alpha }^{\beta } \left( \frac{ {\partial }^{2}f }{
\partial y^{2} } {\phi }^{2} - \frac{ {\partial }^{3}f }{
\partial y^{2} \partial y' } 2 {\phi }^{3} \right) dx
\geq \ \int_{\alpha }^{\beta } \left( \frac{ {\partial }^{2}f }{
\partial y \partial y' } 2 \phi \phi ' \right.$
\item[ ]
$\left. + \ \frac{ {\partial }^{3}f }{ \partial y {\partial
y'}^{2} } 2 {\phi }^{2} \phi ' + \frac{ {\partial }^{2}f }{
{\partial y'}^{2} } \phi \phi '' + \frac{ {\partial }^{3}f }{
\partial y {\partial y'}^{2} } \phi ' {\phi }^{2} + \frac{
{\partial }^{3}f }{ {\partial y'}^{3} } \phi {\phi '}^{2} \right)
dx$, \hfill (18)
\end{enumerate}
so either the above are variational inequalities of motion or the
Lagrangian of motion is not $C^{3}$.
\\ \rm \\
{\bf Proof.} It is enough to notice that $$I''(0) \geq 0. \eqno
(19)$$
\\
{\bf 7. Remark.} Keeping the notation as in the previous
proposition, one may take $$\phi (x) = \lambda \left( (x - \alpha
)(x - \beta ) \right) ^{n}, \eqno (20)$$ where $x \in [ \alpha ,
\beta ]$, while $\lambda > 0$ and $n \in {\bf N} \setminus \{
0,1,2 \} $.
\\ \rm \\
{\bf 8. Example.} If we consider the simple pendulum, where $g$ is
the acce- leration of gravity and $\ell $ is the length of the
weightless thread to the one end of which is connected a particle
of mass $m$, while $\theta $ is the angular displacement as a
function of time $t$, then the Lagrangian of the simple pendulum
is $$L = \frac{1}{2} m {\ell }^{2} {\dot{\theta }}^{2} + m g \ell
\cos \theta \eqno (21)$$ and it is $C^{ \infty }$, so $$\frac{
\partial L}{ \partial t} = 0, \eqno (22)$$
$$\frac{ \partial L}{ \partial \theta } = - mg\ell \sin \theta
\eqno (23)$$ and $$\frac{ \partial L}{ \partial \dot{\theta } } =
m {\ell }^{2} \dot{\theta }, \eqno (24)$$ which imply that the
Euler-Lagrange equation $$\frac{ \partial L}{ \partial \theta } -
\frac{d}{dt} \left( \frac{ \partial L}{ \partial \dot{\theta } }
\right) = 0 \eqno (25)$$ in Proposition 3 takes the form
$$\ddot{\theta } + \frac{g}{\ell } \sin \theta = 0. \eqno (26)$$
So, apart from solving equation (26), $\theta $ must satisfy the
conclusion of Proposition 6 and Remark 7. Since $$\frac{ {\partial
}^{2}L }{ \partial {\theta }^{2} } = - mg \ell \cos \theta , \eqno
(27)$$ $$\frac{ {\partial }^{2}L }{
\partial \theta \partial \dot{\theta } } = 0 \eqno (28)$$ and $$\frac{
{\partial }^{2}L }{ \partial {\dot{\theta }}^{2} } = m {\ell
}^{2}, \eqno (29)$$ while $$\frac{ {\partial }^{3}L }{ \partial
{\theta }^{2} \partial \dot{\theta } } = 0, \eqno (30)$$ $$\frac{
{\partial }^{3}L }{ \partial \theta \partial {\dot{\theta }}^{2} }
= 0 \eqno (31)$$ and $$\frac{ {\partial }^{3}L }{ \partial
{\dot{\theta }}^{3} } = 0, \eqno (32)$$ it follows that if $$\phi
(t) = \lambda \left( (t - \alpha )(t - \beta ) \right) ^{n}, \eqno
(33)$$ where $t \in [ \alpha , \beta ]$, while $\lambda
> 0$ and $n \in {\bf N} \setminus \{ 0,1,2 \} $, then
\begin{enumerate}
\item[ ]
$\int_{\alpha }^{\beta } \left( - mg\ell \cos \theta \cdot {\phi
}^{2} - 0 \cdot 2 {\phi }^{3} \right) dt$
\item[ ]
$\geq \ \int_{\alpha }^{\beta } \left( 0 \cdot 2 \phi \dot{\phi }
+ 0 \cdot 2 {\phi }^{2} \dot{\phi } + m{\ell }^{2} \cdot \phi
\ddot{\phi } + 0 \cdot \dot{\phi } {\phi }^{2} + 0 \cdot \phi
{\dot{\phi }}^{2} \right) dt$ \hfill (34)
\end{enumerate}
or equivalently $$-mg \ell \int_{\alpha }^{\beta } \cos \theta (t)
{\phi (t)}^{2}dt \geq m{\ell }^{2} \int_{\alpha }^{\beta } \phi
(t) \ddot{\phi }(t) dt \eqno (35)$$ or equivalently $$- g
\int_{\alpha }^{\beta } \cos \theta (t) \phi (t)^{2}dt \geq \ell
\int_{\alpha }^{\beta } \phi (t) \ddot{\phi } (t)dt \eqno (36)$$
or equivalently $$- g \int_{\alpha }^{\beta } \cos \theta (t) \phi
(t)^{2}dt \geq \ell \left( \left[ \phi (t) \dot{\phi }(t) \right]
_{t = \alpha }^{t = \beta } - \int_{\alpha }^{\beta } \dot{\phi
}(t) \dot{\phi } (t)dt \right) \eqno (37)$$ or equivalently $$g
\int_{\alpha }^{\beta } \cos \theta (t) \phi (t)^{2}dt \leq \ell
\int_{\alpha }^{\beta } \dot{\phi } (t)^{2} dt. \eqno (38)$$ A
formula for $\theta $, via (26), can be derived from Section 2.1
on pages 69-80 of [5], so if $$0 < {\theta }_{0} < \pi \eqno
(39)$$ and $${\dot{\theta }}_{0} = 0, \eqno (40)$$ then $$t =
\sqrt{ \frac{\ell }{g} } \ln \left( \frac{ \tan \left( \frac{ \pi
}{4} - \frac{ {\theta }_{0} }{4} \right) }{ \tan \left( \frac{ \pi
}{4} - \frac{ \theta }{4} \right) } \right) \eqno (41)$$ and
consequently $$\theta (t) = \pi - 4 \arctan \left( e^{ - t \sqrt{
\frac{g}{\ell } } } \tan \left( \frac{ \pi }{4} - \frac{ {\theta
}_{0} }{4} \right) \right) \eqno (42)$$ must satisfy (38) for all
$\phi $ in question.

\end{document}